\documentclass{ws-procs9x6}  

\begin{document}

\title{Markov Property of Monotone L\'evy Processes}

\author{Uwe FRANZ}

\address{Institut f{\"u}r Mathematik und Informatik \\ Ernst-Moritz-Arndt-Univer\-sit{\"a}t Greifs\-wald \\ Jahnstra{\ss}e 15a, D-17487 Greifswald, Germany \\ Email: franz@uni-greifswald.de}

\author{Naofumi MURAKI}

\address{Mathematics Laboratory \\ Iwate Prefectural University \\ Ta\-kizawa, Iwate 020-0193, Japan \\ Email: muraki@iwate-pu.ac.jp}

\maketitle

\abstracts{Monotone L\'evy processes with additive increments are defined and studied. It is shown that these processes have natural Markov structure and their Markov transition semigroups are characterized using the monotone L\'evy-Khintchine formula.\cite{muraki00} Monotone L\'evy processes turn out to be related to classical L\'evy processes via Attal's ``remarkable transformation.'' A monotone analogue of the family of exponential martingales associated to a classical L\'evy process is also defined.}

\section{Introduction}

One of the remarkable features of quantum probability is the existence of several different notions of independence. The most prominent examples are certainly tensor independence and free independence. Tensor independence is based on the tensor product of algebras and linear functionals. It generalizes the notion of independence used in classical probability and was used to develop the quantum stochastic calculus on the symmetric Fock space, cf.\ Ref.~\refcite{hudson+parthasarathy84,parthasarathy92}. Free independence can be motivated by the free product of groups. Even though it is incompatible with classical probability - there exist no non-trivial examples of classical random variables that are freely independent - it has been used to develop the so-called free probability theory, which has many parallels to classical probability, cf.\ Ref.~\refcite{voiculescu+dykema+nica92,voiculescu97}. Another independence for which limit theorems, infinite-divisibility, convolution semigroups, etc.\ have been studied is boolean independence, see Ref.~\refcite{speicher+woroudi93}.

Monotone independence\cite{muraki01,muraki00} arose from a non-commutative Brownian motion and an interacting Fock space that was introduced in Ref.~\refcite{muraki96a,muraki97a,lu97b}. Recently, Accardi, Ben Ghorbal, and Obata\cite{accardi+ghorbal+obata03} have discovered that monotone independence also appears in the study of the comb product of graphs. Franz\cite{franz01} has shown that monotone independence leads to an associative product of quantum probability spaces and can be used to develop a theory of quantum L\'evy processes analogous to the ones for tensor, free, and boolean independence in Ref.~\refcite{schuermann95b}. The monotone product in not commutative, unlike the tensor, free, or boolean product.

Muraki\cite{muraki02a,muraki02} has proven that there exist only five universal notions of independence in quantum probability. These are tensor, free, boolean, monotone, and anti-monotone independence. Franz\cite{franz03b} subsequently found a construction that reduces boolean, monotone, and anti-monotone independence to tensor independence. As an application the theories of quantum L\'evy processes with boolean, monotone, and anti-monotone increments can be reduced to the theory of L\'evy processes on involutive bialgebras.\cite{schuermann93} This actually implies that for any quantum L\'evy processes with monotone increments their exists a natural family of conditional expectations with respect to which it is Markovian, cf.\ Corollary 4.4 in Ref.~\refcite{franz03b}. But the construction is algebraic and does not give a direct expression for the semigroup of Markovian transition kernel.

In the present article we give a description of the semigroups of Markovian transition operators and kernels of quantum L\'evy processes with additive monotonically independent increments. This is similar to the program carried out by Biane\cite{biane98} for quantum L\'evy processes with additive or multiplicative freely independent increments. But unlike the free product, the monotone product does not preserve traces, see Remark \ref{mon-non-trace-non-unit}. Therefore the existence of the conditional expectations does not follow from the theory of von Neumann algebras, but relies on the explicit construction given in Ref.~\refcite{franz03b}. It turns out that monotone L\'evy processes can be obtained from the Fock space realization of classical L\'evy processes via Attal's ``remarkable transformation,''\cite{attal01} see the proof of Theorem \ref{mon-theo-mon-levy}. The Fock space realization of a classical L\'evy process defined in Equation \eqref{mon-eq-class-levy} is a quantum semi-martingale in the sense of Ref.~\refcite{attal01} only if the L\'evy measure has compact support. In this case it turns out that the associated monotone L\'evy process is even bounded for all $t\ge 0$. We also associate a family of martingales to a monotone L\'evy process that are analogous to the exponential martingales of a classical L\'evy process, see Theorem \ref{mon-theo-martingale}.

\section{Monotone Independence}

In this section we present the definition of monotone independence and its main properties. 

\begin{definition}\label{mon-def-mon-indep}
Let $\mathcal{H}$ be a Hilbert space, $\Omega\in\mathcal{H}$ a unit vector, and define a state $\Phi:\mathcal{B}(\mathcal{H})\to\mathbb{C}$ on the algebra of bounded operators on $\mathcal{H}$ by
\[
\Phi(X)=\langle \Omega, X\Omega\rangle, \qquad\mbox{ for } X\in\mathcal{B}(\mathcal{H}).
\]
Subalgebras $\mathcal{A}_1,\ldots,\mathcal{A}_r\subseteq\mathcal{B}(\mathcal{H})$ are called {\em monotonically independent} w.r.t.\ $\Phi$, if the following two conditions are satisfied.
\begin{itemize}
\item[(a)]
For all $X\in\mathcal{A}_i$, $Y\in\mathcal{A}_j$, $Z\in\mathcal{A}_k$ with $j>\max(i,k)$, we have
\[
XYZ = \Phi(Y)XZ.
\]
\item[(b)]
For all $X_1\in\mathcal{A}_{i_1}$, $\ldots$ , $X_n\in\mathcal{A}_{i_n}$, $Y\in\mathcal{A}_j$, $Z_1\in\mathcal{A}_{k_1}$, $\ldots$, $Z_m\in\mathcal{A}_{k_m}$ with $i_1>i_2>\cdots >i_n$, $k_1>k_2>\ldots>k_m$, and $j<\min(i_n,k_m)$, we have
\[
\Phi(X_1\cdots X_nYZ_m\cdots Z_1)=\Phi(X_1)\cdots\Phi(X_n)\Phi(Y)\Phi(Z_m)\cdots\Phi(Z_1).
\] 
\end{itemize}
Operators $X_1,\ldots,X_r\in\mathcal{B}(\mathcal{H})$ are called monotonically independent w.r.t.\ $\Omega$, if the subalgebras $\mathcal{A}_i={\rm alg}(X_i)={\rm span}\{X_i^k|k=1,2,\ldots\}$, $i=1,\ldots,r$, they generate are monotonically independent.
\end{definition}

\begin{remark}\label{mon-non-unit}
We do not require the subalgebras to be unital. If the $i^{\rm th}$ algebra $\mathcal{A}_i$ contains the unit, then $\Phi|_{\mathcal{A}_k}$ is a homomorphism for $k>i$, since
\[
\Phi(X_1X_2)=\Phi(X_1\mathbf{1}X_2)=\Phi(\mathbf{1})\Phi(X_1)\Phi(X_2)=\Phi(X_1)\Phi(X_2)
\]
for all $X_1,X_2\in\mathcal{A}_k$.
\end{remark}

We will call a triple $(\mathcal{A},\mathcal{H},\Omega)$ consisting of a Hilbert space $\mathcal{H}$, a unit vector $\Omega\in\mathcal{H}$, and a subalgebra $\mathcal{A}\subseteq\mathcal{B}(\mathcal{H})$ a {\em quantum probability space}. By an {\em operator process} we shall mean an indexed family $(X_t)_{t\in I}$ of elements of some quantum probability space. A {\em quantum random variable} is a homomorphism $j:\mathcal{B}\to \mathcal{A}$ from some $*$-algebra into a quantum probability space and a {\em quantum stochastic process} is an indexed family $(j_t)_{t\in I}$ of quantum random variables.

\begin{proposition}\label{mon-prop-mon-prod}
Let $(\mathcal{A}_i,\mathcal{H}_i,\Omega_i)$, $i=1,2$, be two quantum probability spaces, and denote the states associated to $\Omega_1$ and $\Omega_2$ by $\Phi_1$ and $\Phi_2$, respectively.

Then there exists a quantum probability space $(\mathcal{A},\mathcal{H},\Omega)$ and two injective state-preserving homomorphisms $J_i:\mathcal{A}_i\to\mathcal{A}$, $i=1,2$, such that the images $J_1(\mathcal{A}_1)$ and $J_2(\mathcal{A}_2)$ are monotonically independent w.r.t.\ $\Omega$.
\end{proposition}
\begin{proof}
We set $\mathcal{H}=\mathcal{H}_1\otimes\mathcal{H}_2$ and $\Omega=\Omega_1\otimes\Omega_2$. Denote by $P_2$ the orthogonal projection on $\mathbb{C}\Omega_2\subseteq\mathcal{H}_2$.

We define the embeddings $J_i:\mathcal{A}_i\to\mathcal{B}(\mathcal{H})$ by
\begin{eqnarray*}
J_1(X) &=& X\otimes P_2, \qquad \mbox{ for } X\in\mathcal{A}_1, \\
J_2(X) &=& \mathbf{1}\otimes X, \qquad \mbox{ for } X\in\mathcal{A}_2.
\end{eqnarray*}
For $\mathcal{A}$ we take the subalgebra generated by $J_1(\mathcal{A}_1)$ and $J_2(\mathcal{A}_2)$. It is clear that $J_1$ and $J_2$ are injective, state-preserving homomorphism.

A simple calculation shows that $J_1(\mathcal{A}_1)$ and $J_2(\mathcal{A}_2)$ are monotonically independent w.r.t.\ $\Omega$. E.g., for products of the form $J_1(X_1)J_2(Y)J_1(X_2)$, $X_1,X_2\in\mathcal{A}_1$, $Y\in\mathcal{A}_2$, we get
\begin{eqnarray*}
J_1(X_1)J_2(Y)J_1(X_2) &=& (X_1\otimes P_2)(\mathbf{1}\otimes Y)(X_1\otimes P_2) = (X_1X_2)\otimes P_2YP_2 \\
&=& \Phi\big(J_2(Y)\big) J_1(X_1)J_1(X_2).
\end{eqnarray*}
On the other hand, for $J_2(Y_1)J_1(X)J_2(Y_2)$, $X\in\mathcal{A}_1$, $Y_1,Y_2\in\mathcal{A}_2$, we get
\begin{eqnarray*}
\Phi\big(J_2(Y_1)J_1(X)J_2(Y_2)\big) &=& \langle\Omega_1\otimes\Omega_2,(\mathbf{1}\otimes Y_1)(X\otimes P_2)(\mathbf{1}\otimes Y_2)\Omega_1\otimes\Omega_2\rangle \\
&=& \langle\Omega_1\otimes\Omega_2,X\otimes (Y_1PY_2)\Omega_1\otimes\Omega_2\rangle \\
&=& \Phi_1(X)\Phi_2(Y_1)\Phi_2(Y_2) \\
&=& \Phi\big(J_2(Y_1)\big)\Phi\big(J_1(X)\big)\Phi\big(J_2(Y_2)\big).
\end{eqnarray*}
\end{proof}

We will call the quantum probability space  $(\mathcal{A},\mathcal{H},\Omega)$ constructed in the previous proposition the {\em monotone product} of $(\mathcal{A}_1,\mathcal{H}_1,\Omega_1)$ and $(\mathcal{A}_2,\mathcal{H}_2,\Omega_2)$. When there is no danger of confusion, we shall identify the algebras $\mathcal{A}_1$ and $\mathcal{A}_2$ with their images $J_1(\mathcal{A}_1)$ and $J_2(\mathcal{A}_2)$, respectively.

\begin{remark}\label{mon-non-trace-non-unit}
\begin{itemize}
\item[(a)]
The monotone product is associative and can be extended to more than two factors, see also Ref.~\refcite{franz01}. But it is not commutative.
\item[(b)]
The state $\Phi$ on $\mathcal{A}={\rm alg}(\mathcal{A}_1,\mathcal{A}_2)$ depends only on $\Phi_1|_{\mathcal{A}_1}$ and $\Phi_2|_{\mathcal{A}_2}$.
\item[(c)]
The embedding $J_1:\mathcal{A}_1\to\mathcal{A}$ is not unital. It is a consequence of Remark \ref{mon-non-unit} that it is impossible to get a unital embedding of the first algebra, if the state on the second algebra is not a homomorphism.
\item[(d)]
The product is not trace-preserving. If $\Phi_1|_{\mathcal{A}_1}$ is not identically equal to zero, then the calculation
\[
\Phi_1(X)\Phi_2(Y_1Y_2)=\Phi(XY_1Y_2)=\Phi(Y_2XY_1)=\Phi_1(X)\Phi_2(Y_1)\Phi_2(Y_2)
\]
for all $X\in\mathcal{A}_1$, $Y_1,Y_2\in\mathcal{A}_2$ shows that $\Phi$ can only be a trace on $\mathcal{A}$, if $\Phi_2|_{\mathcal{A}_2}$ is a homomorphism. Note that we identified the elements of $\mathcal{A}_1$ and $\mathcal{A}_2$ with their images under the embedding $J_1$ and $J_2$ to simplify the notation.
\end{itemize}
\end{remark}
We now recall several results from Ref.~\refcite{muraki00}, which will be important for the following sections.

By $\mathbb{C}^+=\{z\in\mathbb{C}|\Im z>0\}$ we denote the upper half plane.

Denote by $R_X(z)=(z-X)^{-1}$ the resolvent of an operator $X$. Let $X\in\mathcal{B}(\mathcal{H})$ now be self-adjoint. We will need the {\em reciprocal Cauchy transform}
\[
H_X(z) = \frac{1}{\Phi\big((z-X)^{-1}\big)}, \qquad z\in\mathbb{C}^+,
\]
of the spectral measure of $X$ evaluated in the state $\Phi$. We use the same notation for the reciprocal Cauchy transform
\[
H_\nu(z)=\frac{1}{\int_\mathbb{R}\frac{1}{z-x}{\rm d}\nu(x)}, \qquad z\in\mathbb{C}^+,
\]
of a probability measure $\nu$ on $\mathbb{R}$.

\begin{theorem} (Theorem 3.1 in Ref.~\refcite{muraki00})
Let $X_1,\ldots,X_n\in\mathcal{B}(\mathcal{H})$ be monotonically independent self-adjoint operators. Then
\[
H_{X_1+\cdots+X_n}(z) = H_{X_1}\big(\cdots \big(H_{X_n}(z)\big)\cdots\big)
\]
for $z\in\mathbb{C}^+$.
\end{theorem}

This theorem suggests to define the monotone convolution of probability measures as follows.
\begin{definition}(Definition 3.2 and Theorem 3.5 in Ref.~\refcite{muraki00})
For any pair of probability measures $\mu,\nu$ on $\mathbb{R}$ there exists a unique probability measure $\lambda$ on $\mathbb{R}$, whose reciprocal Cauchy transform is given by
\[
H_\lambda = H_\mu\big(H_\nu(z)\big),\qquad z\in\mathbb{C}^+.
\]
We will call $\lambda$ the {\em monotone convolution} of $\mu$ and $\nu$ and denote it by $\lambda=\mu\triangleright\nu$.
\end{definition}

In the proof of Theorem 3.5 in Ref.~\refcite{muraki00}, Muraki gives the following formula for the monotone convolution,
\begin{equation}\label{mon-formula-mon-conv}
\mu\triangleright\nu (\cdot)= \int_{\mathbb{R}} \nu_y(\cdot){\rm d}\mu(y)
\end{equation}
where the measures $\mu_y$, $y\in\mathbb{R}$, are given by their reciprocal Cauchy transforms
\[
H_{\nu_y}(z)=H_\nu(z)-y, \qquad z\in\mathbb{C}^+.
\]
The measures $\nu_y$, $y\in\mathbb{R}$, can also be defined by $\nu_y=\delta_y\triangleright \nu$. The monotone convolution is associative, affine in the first argument, and weakly continuous in both arguments. Note that it is not commutative.

A {\em monotone convolution semigroup} is a weakly continuous family $(\mu_t)_{t\ge 0}$ of probability measures on $\mathbb{R}$ such that $\mu_{s+t}=\mu_s\triangleright\mu_t$ for all $s,t\ge 0$ and $\mu_0=\delta_0$. Monotone convolution semigroups are in one-to-one correspondence with continuous families $(H_t)_{t\ge0}$ of reciprocal Cauchy transforms that form a semigroup w.r.t.\ composition. These semigroups are even differentiable with respect to $t$ and the derivative
\[
A(z)=\lim_{t\searrow 0} \frac{1}{t}\big(H_t(z)-z\big),\qquad z\in\mathbb{C}^+
\]
is called their {\em generator}. Muraki has classified all monotone convolution semigroups in terms of their generators, see Theorem 4.5 in Ref.~\refcite{muraki00}. Since we will only consider monotone convolution semigroups of probability measures with compact support in this paper, we will need Theorem 5.1 from Ref.~\refcite{muraki00}.

\begin{theorem}{\bf (L\'evy-Khintchine formula for the monotone convolution)} \label{mon-theo-levy-khintchine}
Let $(\mu_t)_{t\ge 0}$ be a weakly continuous family of probability measures and assume $\mu_t\not=\delta_0$ for some $t>0$. Then the following are equivalent.
\begin{itemize}
\item[(a)]
$(\mu_t)_{t\ge0}$ is a monotone convolution semigroup of compactly supported measures on $\mathbb{R}$.
\item[(b)]
There exists a unique pair $(a,\rho)\not=(0,0)$ of a real number $a$ and a compactly supported finite measure $\rho$ on $\mathbb{R}$ such that the Pick function
\[
A(z)=a+\int_\mathbb{R}\frac{1}{x-z}{\rm d}\rho(x)
\]
generates $H_{\mu_t}$ in the sense that for $z\in\mathbb{C}^+$, $H_{\mu_t}(z)=w$ is the unique solution $w\in\mathbb{C}^+$ of the equation
\[
\int_z^w\frac{{\rm d}\zeta}{A(\zeta)} = t.
\]
\end{itemize}
\end{theorem}

\begin{remark}
Composition semigroups of holomorphic functions on a half plane or on the disc and their generators were also studied in Ref.~\refcite{berkson+porta78,goryainov93}.
\end{remark}

We will call the pair $(a,\rho)$ the {\em characteristic pair} of the monotone convolution semigroup $(\mu_t)_{t\ge 0}$, $-a$ its {\em drift coefficient}, $\rho(\{0\})$ its {\em diffusion coefficient},  and $\frac{1}{x^2}\big(\rho-\rho(\{0\})\delta_0\big)$ its {\em L\'evy measure}.

\section{Conditional Expectations}

We will first introduce ``conditional expectations'' for the monotone product of two quantum probability spaces.

\begin{proposition}\label{mon-prop-con-exp}
Let $(\mathcal{A}_i,\mathcal{H}_i,\Omega_i)$, $i=1,2$, be two quantum probability spaces and denote by $(\mathcal{A},\mathcal{H},\Omega)$ their monotone product. Then there exists a linear map $E_1:\mathcal{A}\to\mathcal{A}_1$ such that
\begin{itemize}
\item[(a)]
$E_1\big(J_1(X)YJ_1(Z)\big)=XE_1(Y)Z$ for all $X,Z\in\mathcal{A}_1$, $Y\in\mathcal{A}$,
\item[(b)]
$\Phi_1\circ E_1 =\Phi$,
\item[(c)]
$E_1\circ J_1 = {\rm id}$,
\item[(d)]
$E_1$ is completely positive,
\item[(e)]
if $\mathcal{A}_1$ and $\mathcal{A}$ contain units $\mathbf{1}_1$ and $\mathbf{1}$, then we have $E_1(\mathbf{1})=\mathbf{1}_1$.
\end{itemize}
We will call $E_1$ the {\em conditional expectation} from $(\mathcal{A},\mathcal{H},\Omega)$ onto $(\mathcal{A}_1,\mathcal{H}_1,\Omega_1)$.
\end{proposition}
\begin{proof}
Let $P$ denote the orthogonal projection onto $\mathcal{H}_1\otimes\Omega_2\subseteq\mathcal{H}$, i.e.\ $P=\mathbf{1}\otimes P_2$. Since
\[
P (X\otimes Y) P = X \otimes P_2YP_2 = \Phi(Y) X\otimes P_2,
\]
for $X\in\mathcal{B}(\mathcal{H}_1)$, $Y\in\mathcal{B}(\mathcal{H}_2)$, we see that the image of the map from $\mathcal{A}$ to $\mathcal{B}(\mathcal{H})$ defined by $Z\mapsto PZP$ lies in $J_1(\mathcal{A}_1)$. Therefore we can define $E_1:\mathcal{A}\to\mathcal{A}_1$ by $E_1(Z)=J_1^{-1}(PZP)$.

It is straight-forward to verify that $E_1$ satisfies all the properties listed in the proposition.
\end{proof}

\begin{remark}
The map $Z\mapsto PZP$ from $\mathcal{B}(\mathcal{H})$ to $\mathcal{B}(\mathcal{H}_1)\otimes P_2$ is basically the conditional expectation for the tensor product  of quantum probability spaces. Note that it is impossible to define a conditional expectation onto $\mathcal{A}_2$ in the same way, since the image of the map $Z\mapsto P'ZP'$, where $P'$ is the orthogonal expectation onto $\Omega_1\otimes \mathcal{H}_2\subseteq\mathcal{H}$, does not lie inside $J_2(\mathcal{A}_2)$.
\end{remark}

\begin{proposition}\label{mon-prop-cond-exp}
Let $X_1$ and $X_2$ be two self-adjoint operators on $\mathcal{H}_i$ and fix unit vectors $\Omega_i\in\mathcal{H}_i$, $i=1,2$. Denote by $\mathcal{A}_i$ the $C^*$-algebras $\mathcal{A}_i=\{f(X_i)|f\in C_b(\mathbb{R})\}$ generated by $X_1$ and $X_2$, respectively. Let $(\mathcal{A},\mathcal{H},\Omega)$ be the monotone product of $(\mathcal{A}_i,\mathcal{H}_i,\Omega_i)$, $i=1,2$, and denote the images of $X_1$ and $X_2$ in $\mathcal{A}$ again by $X_1$ and $X_2$.

Then we have
\[
E_1\left(\frac{1}{z-(X_1+X_2)}\right) = \frac{1}{H_{X_2}(z)-X_1}
\]
for $z\in\mathbb{C}^+$.
\end{proposition}
\begin{proof}
The idea of the proof is the same as for Proposition 3.2 in Ref.~\refcite{biane98}. Denote again by
\[
R_X(z)=(z-X)^{-1}
\]
the resolvent of an operator $X$. Then we have
\[
R_{X_1+X_2}(z)=R_{X_2}(z)\big(1-X_1R_{X_2}(z)\big)^{-1}
\]
As in Ref.~\refcite{biane98} we can expand this expression into a norm convergent series
\[
R_{X_1+X_2}(z) = \sum_{k=0}^\infty R_{X_2}(z)\big(X_1R_{X_2}(z)\big)^k
\]
for $|z|>||X_1||+||X_2||$. Therefore
\begin{eqnarray*}
PR_{X_1+X_2}(z)P &=& \sum_{k=0}^\infty \Big(\Phi\big(R_{X_2}(z)\big)\Big)^{k+1}X_1^kP \\
&=& \sum_{k=0}^\infty \big(H_{X_2}(z)\big)^{-k-1}X_1^k P  = \frac{P}{H_{X_2}(z)-X_1},
\end{eqnarray*}
and
\[
E_1\big(R_{X_1+X_2}(z)\big) =\frac{1}{H_{X_2}(z)-X_1}.
\]

By uniqueness of analytic continuation follows that the identity holds for all $z\in \mathbb{C}\backslash[-||X_1||-||X_2||,||X_1||+||X_2||]$.
\end{proof}
\begin{corollary} \label{mon-cor-cond-exp}
For all $f\in C_b(\mathbb{R})$, we have
\[
E_1\big(f(X_1+X_2)\big) = (Tf)(X_1),
\]
where the operator $T$ is defined by
\[
Tf(x)=\int_\mathbb{R}f(y)\mu_{2,x}({\rm d}y), \qquad x\in\mathbb{R},
\]
with the measures $\mu_{2,x}$ determined by their reciprocal Cauchy transforms,
\[
H_{\mu_{2,x}}(z)= H_{X_2}(z)-x, \qquad \mbox{ for } z\in\mathbb{C}^+.
\] 
\end{corollary}
\begin{remark}
We can also prove this formula using only the fact that $T$ does not depend on the distribution of $X_1$ and Formula \eqref{mon-formula-mon-conv}.

Denote the conditional expectation of $f(X_1+X_2)$ by $Tf(X_1)=E_1\big(f(X_1+X_2)\big)$. Since the conditional expectation preserves expectations, we
get
\[
\Phi\big(Tf(X_1)\big)=\Phi\big(f(X_1+X_2)\big).
\]
Denoting by $\mu_1$ and $\mu_2$ the distributions of $X_1$ and $X_2$ w.r.t.\  $\Omega$, this becomes
\[
\int_\mathbb{R} Tf(x)\mu_1({\rm d}x) = \int_{\mathbb{R}} f(x)\mu_1\triangleright\mu_2({\rm d}x).
\]
Substituting $\mu_1\triangleright\mu_2$ with Formula \eqref{mon-formula-mon-conv}, we obtain
\[
\int_\mathbb{R} Tf(x)\mu_1({\rm d}x) = \int_\mathbb{R} f(x) \,{\rm d}\!\!\int_\mathbb{R} \mu_{2,y}(x) \mu_1({\rm d}y).
\]
Since this has to hold for all probability measures $\mu_1$, we get the desired result.
\end{remark}

\section{Monotone L\'evy Processes}

\begin{definition}\label{mon-def-mon-levy}
A family of self-adjoint operators $(X_t)_{t\ge0}\subseteq\mathcal{B}(\mathcal{H})$ is called {\em monotone L\'evy process} w.r.t.\ $\Omega\in\mathcal{H}$, $||\Omega||=1$, if the following conditions are satisfied.
\begin{itemize}
\item[(a)]
(Independence of increments) For all $n\in\mathbb{N}$ and $0\le t_1\le \cdots\le t_n$, the operators
\[
X_{t_1}, X_{t_2}-X_{t_1},\ldots,X_{t_n}-X_{t_{n-1}},
\]
are monotonically independent w.r.t.\ $\Omega$.
\item[(b)]
(Stationarity)
The distribution of an increment $X_t-X_s$ w.r.t.\ to the state $\Phi(\cdot)=\langle\Omega,\cdot\,\Omega\rangle$ depends only on $t-s$, i.e.\
\[
\Phi\big(f(X_t-X_s)\big) = \Phi\big(f(X_{t+h}-X_{s+h})\big)
\]
for all $0\le s\le t$, $h\ge 0$, and $f\in C_b(\mathbb{R})$.
\item[(c)]
(Weak continuity) $\lim_{t\searrow 0} \Phi\big(f(X_t)\big)=f(0)$ for all $f\in C_b(\mathbb{R})$.
\end{itemize}
\end{definition}

We will call two monotone L\'evy processes $(X_t)_{t\ge 0}$ and $(X'_t)_{t\ge 0}$, defined on $(\mathcal{H},\Omega)$ and $(\mathcal{H}',\Omega')$, {\em equivalent} if their marginal distributions coincide, i.e.\
\[
\langle \Omega, f(X_t-X_s)\Omega\rangle = \langle \Omega', f(X'_t-X'_s)\Omega'\rangle
\]
for all $0\le s\le t$, $f\in C_b(\mathbb{R})$. Due to the monotone independence of the increments this implies that all finite joint distributions also coincide.

Denote by $\mathbb{C}[x]$ the algebra of polynomials generated by one self-adjoint indeterminate $x=x^*$. It is a dual group with the comultiplication $\Delta:\mathbb{C}[x]\to \mathbb{C}[x]\coprod\mathbb{C}[x]\cong\mathbb{C}[x',x'']$ defined by $\Delta(x)=x'+x''$, see Ref.~\refcite{franz01,franz03b}.

The conditions in Definition \ref{mon-def-mon-levy} imply that the homomorphisms $j_{st}:\mathbb{C}[x]\to\mathcal{B}(\mathcal{H})$ defined by
\[
j_{st}(\mathbf{1})={\rm id},\qquad j_{st}(x^k)=(X_t-X_s)^k
\]
for $0\le s\le t$ and $k=1,2,\ldots$ form a monotone L\'evy process on $\mathbb{C}[x]$ in the sense of Definition 2.5 in Ref.~\refcite{franz03b}. In Ref.~\refcite{franz01,franz03b} it was shown that such a processes is uniquely characterized by its convolution semigroup $(\varphi_t)_{t\ge 0}$ of states on $\mathbb{C}[x]$ given by
\[
\varphi_t(x^k)=\langle \Omega, X^k_t\Omega\rangle, \qquad t\ge0, \quad k\ge 0.
\]

Conversely, a monotone L\'evy processes $(j_{st})_{0\le s\le t}$ on $\mathbb{C}[x]$ in the sense of Definition 2.5 in Ref.~\refcite{franz03b} defines a monotone L\'evy process in the sense of Definition \ref{mon-def-mon-levy} by
\[
X_t=j_{0t}(x),\qquad t\ge 0,
\]
if these operators are bounded (and hence self-adjoint) for all $t\ge 0$. 

Since we are only considering bounded operators, the marginal distributions $(\mu_t)_{t\ge 0}$ defined by
\[
\int f(x){\rm d}\mu_t(x) = \langle\Omega, f(X_t)\Omega\rangle, \qquad \mbox{ for all } f\in C_b(\mathbb{R}),
\]
for $t\ge 0$ have compact support contained in $[-||X_t||,||X_t||]$ and are uniquely determined by their moments. Therefore we obtain the following classification of monotone L\'evy processes.

\begin{proposition}\label{mon-prop-one-to-one}
We have a one-to-one correspondence between monotone L\'evy processes $(X_t)_{t\ge 0}$ (modulo equivalence) and monotone convolution semigroups $(\mu_t)_{t\ge 0}$ of compactly supported probability measures.
\end{proposition}

We will now apply the results of Ref.~\refcite{franz01,franz03b} to show how the monotone L\'evy process associated to a monotone convolution semigroup of compactly supported measures with characteristic pair $(a,\rho)$ can be constructed.

Let $\mu$ be a compactly supported probability measure on $\mathbb{R}$, ${\rm supp}\, \mu\subseteq[-M,M]$, with $M>0$, and define $\varphi_\mu:\mathbb{C}[x]\to\mathbb{C}$ by
\[
\varphi(P)=\int_\mathbb{R}P(x){\rm d}\mu(x), \qquad P\in\mathbb{C}[x].
\]
Then $\varphi_\mu$ is uniquely determined by the generating function
\[
\sum_{k=0}^\infty \varphi(x^k) z^{-k-1} = \int_{\mathbb{R}}\frac{1}{z-x}{\rm d}\mu(x) = \frac{1}{H_\mu(z)},
\]
 for $|z|>M$. It follows that the generator
\[
L=\left.\frac{{\rm d}}{{\rm d}t}\right|_{t=0}\varphi_t
\]
of the convolution semigroup of states $(\varphi_t)_{t\ge0}$, $\varphi_t=\varphi_{\mu_t}$ for $t\ge 0$, associated to a monotone convolution semigroup $(\mu_t)_{t\ge 0}$, can be characterized by the generating function
\[
\sum_{k=0}^\infty L(x^k)z^{-k-1} = -\frac{A(z)}{z^2} =-\frac{a}{z^2}+\frac{1}{z^2}\int_\mathbb{R}\frac{1}{z-x}{\rm d}\rho(x)
\]
for sufficiently large $|z|$. Therefore we get
\[
L(x^k)=\left\{\begin{array}{lcl}
0 & \mbox{ if } & k=0, \\
-a & \mbox{ if } & k=1, \\
\int_{\mathbb{R}} x^{k-2}{\rm d}\rho(x) & \mbox{ if } & k=2,3,\ldots.
\end{array}\right.
\]
Note that $L:\mathbb{C}[x]\to\mathbb{C}$ extends to a unique continuous functional on $C_b^2(\mathbb{R})$, which can be given by
\begin{equation}\label{mon-gen-L}
Lf=-af'(0)+\int_{\mathbb{R}}\big(f(x)-f(0)-xf'(0)\big)\frac{{\rm d}\rho(x)}{x^2}
\end{equation}
for $f\in C_b^2(\mathbb{R})$.

Recall that a {\em Sch\"urmann triple} $(\pi,\eta,L)$ on a $*$-algebra $\mathcal{B}$ with a character $\varepsilon:\mathcal{B}\to\mathbb{C}$ over some (pre-) Hilbert space $H$ consists of a 
\begin{itemize}
\item
a $*$-representation $\pi$ of $\mathcal{B}$ on $H$,
\item
a $\pi$-$\varepsilon$-cocycle $\eta$, i.e.\ a linear map $\eta:\mathcal{B}\to H$ such that
\begin{equation}\label{mon-cocycle-eta}
\eta(ab)=\pi(a)\eta(b)+\eta(a)\varepsilon(b)
\end{equation}
for all $a,b\in\mathcal{B}$, and 
\item
a hermitian linear functional $L:\mathcal{B}\to\mathbb{C}$ such that
\begin{equation}\label{mon-coboundary-L}
L(ab)=\varepsilon(a)L(b)+\langle\eta(a^*),\eta(b)\rangle+L(a)\varepsilon(b)
\end{equation}
holds for all $a,b\in\mathcal{B}$.
\end{itemize}

\begin{lemma}\label{mon-lemma-triple}
Let $H=L^2(\mathbb{R},\rho)$ and set $\pi(f)=M_f$ (the multiplication operator on $L^2(\mathbb{R},\rho)$, $M_f(g)=fg$) and
\[
\eta(f)=\left\{\begin{array}{lcl}
\frac{f(x)-f(0)}{x} & \mbox{ if } & x\not=0, \\
f'(0) & \mbox{ if } & x=0,
\end{array}\right.
\]
for $f\in C^2_b(\mathbb{R})$.

Then $(\pi,\eta,L)$ with $L$ as in Equation \eqref{mon-gen-L} defines a Sch\"urmann triple on $C^2_b(\mathbb{R})$ with the character $\varepsilon(f)=f(0)$ for $f\in C_b^2(\mathbb{R})$.
\end{lemma}
\begin{proof}
Clearly, $\pi$ is a representation.

Let $f,g\in C_b^2(\mathbb{R})$, then we have
\begin{eqnarray*}
\eta(fg) &=&\left\{\begin{array}{lcl}
\frac{(fg)(x)-(fg)(0)}{x} =f(x)\frac{g(x)-g(0)}{x}+g(0)\frac{f(x)-f(0)}{x} & \mbox{ if } & x\not=0 \\
(fg)'(0) = f'(0)g(0)+f(0)g'(0) & \mbox{ if } & x=0
\end{array}\right. \\
&=& \pi(f)\eta(g)+\eta(f)\varepsilon(g)
\end{eqnarray*}
i.e.\ Equation \eqref{mon-cocycle-eta} is satisfied. Furthermore, we get
\begin{eqnarray*}
L(fg) &=& -a (fg)'(0) + \int_\mathbb{R} \big((fg)(x)-(fg)(0)-x(fg)'(0)\big))\frac{{\rm d}\rho(x)}{x^2} \\
&=& f(0)\left(-ag'(0) + \int_\mathbb{R} \big(g(x)-g(0)-xg'(0)\big)\frac{{\rm d}\rho(x)}{x^2}\right) \\
&& + \int_\mathbb{R}\frac{\big(f(x)-f(0)\big)\big(g(x)-g(0)\big)}{x^2}{\rm d}\rho(x) \\
&& + g(0) \left(-af'(0) + \int_\mathbb{R} \big(f(x)-f(0)-xf'(0)\big)\frac{{\rm d}\rho(x)}{x^2}\right) \\
&=& \varepsilon(f)L(g)+\langle\eta(\overline{f}),\eta(g)\rangle+L(f)\varepsilon(g),
\end{eqnarray*}
and therefore Equation \eqref{mon-coboundary-L} is also satisfied.
\end{proof}

For $0\le s\le t$, let $P_{]s,t[^c}:L^2(\mathbb{R}_+,H)\to L^2(\mathbb{R}_+,H)$ be the orthogonal projection onto $L^2\big(]s,t[^c,H)\subseteq L^2(\mathbb{R}_+,H)$, i.e.\ $P_{]s,t[^c}f=\mathbf{1}_{]s,t[^c}f$. We denote by $\Gamma_{]s,t[^c}=\Gamma(P_{]s,t[^c})\in\mathcal{B}\Big(\Gamma\big(L^2(\mathbb{R}_+,H)\big)\Big)$ , $0\le s\le t$, the second quantization of $P_{]s,t[^c}$. For the case $t=\infty$, we introduce the shorter notation $P_{s]}=P_{[s,\infty[^c}$ and $\Gamma_{s]}=\Gamma_{]s,\infty[^c}$.

\begin{theorem}\label{mon-theo-mon-levy}
Let $a$ be a real number, $\rho$ a compactly supported finite measure on $\mathbb{R}$, and $H=L^2(\mathbb{R},\rho)$. Denote by $M_x$ the operator $H\ni f\mapsto xf\in H$ and by $\mathbf{1}_\mathbb{R}\in H$ the constant function with value one. Then the quantum stochastic differential Equation
\begin{equation}\label{mon-bialg-qsde}
\tilde{X}_{st}= \int_s^t\left( -\tilde{X}_{s\tau}{\rm d}\Lambda_\tau({\rm id})+{\rm d}\Lambda_\tau(M_x)+{\rm d}A^+_\tau(\mathbf{1}_\mathbb{R})+ {\rm d}A_\tau(\mathbf{1}_\mathbb{R})-a{\rm d}\tau\right)
\end{equation}
has a unique bounded solution.

Furthermore, $(X_t)_{t\ge0}$ with
\[
X_t=\tilde{X}_{0t}\Gamma_{t]}, \qquad t\ge 0,
\]
is a monotone L\'evy process w.r.t.\ the vacuum vector $\Omega$. The monotone convolution semigroup associated to $(X_t)_{t\ge 0}$ has characteristic pair $(a,\rho)$.

If $(X'_t)_{t\ge 0}$ is another monotone L\'evy processes whose convolution semigroup has characteristic pair $(a,\rho)$, then $(X_t)_{t\ge 0}$ and $(X'_t)_{t\ge0}$ are equivalent.
\end{theorem}
\begin{proof}
Denote by $\mathcal{B}$ involutive algebra generated freely (as an algebra) by two self-adjoint generators $x$ and $p$, with the coalgebra structure defined by
\begin{gather*}
\Delta(x)=x\otimes p + 1\otimes x , \qquad \Delta(p)=p\otimes p, \\
\varepsilon(x)=0, \qquad \varepsilon(p)=1.
\end{gather*}
It follows from Section 3.2 in Ref.~\refcite{franz03b} and Sch\"urmann's representation theorem\cite{schuermann93} that Equation \eqref{mon-bialg-qsde} has a solution on some dense invariant subspace of $\Gamma\big(L^2(\mathbb{R}_+,H)\big)$ and that
\[
\tilde{\jmath}_{st}(x)=\tilde{X}_{st} \quad\mbox{ and } \quad \tilde{\jmath}_{st}(p)=\Gamma_{]s,t[^c}
\]
defines a L\'evy process $(\tilde{\jmath}_{st})_{0\le s\le t}$ on the involutive bialgebra $\mathcal{B}$ w.r.t.\ to the vacuum state in the sense of Ref.~\refcite{schuermann93}.

The operator process $(Y_t)_{t\ge 0}$ defined by
\[
Y_{t}=\Lambda_{t}(M_x)+A^+_{t}(\mathbf{1}_\mathbb{R})+ A_{t}(\mathbf{1}_\mathbb{R})-at{\rm id}, \qquad t\in\mathbb{R}_+,
\]
is a quantum stochastic integral process with bounded coefficients and belongs therefore to the algebra $\mathcal{S}'$ of possibly unbounded quantum semimartingales introduced by Attal,\cite{attal01} see also Ref.~\refcite{belton03}. The operator process $(\tilde{X}_{0t})_{t\ge0}$ satisfies the quantum stochastic differential equation
\begin{equation}\label{mon-eq-class-levy}
\tilde{X}_{0t} = Y_t - \int_0^t \tilde{X}_{0s}{\rm d}\Lambda_s({\rm id}), \qquad t\in\mathbb{R}_+.
\end{equation}
By Proposition 9 from Ref.~\refcite{attal01} and Proposition 37 from Ref.~\refcite{belton03} the solution of this quantum stochastic differential equation is unique and equal to the image $\mathcal{D}Y$ of $(Y_t)_{t\ge 0}$ under Attal's ``remarkable transformation'' $\mathcal{D}$, see also Ref.~\refcite{belton03} for the generalisation to a Fock space with more than one degree of freedom. Therefore $(\tilde{X}_{0t})_{t\ge 0}$ belongs to the algebra $\mathcal{S}$ of bounded quantum semimartingales in the sense of Ref.~\refcite{attal01,belton03}. The operators $\tilde{X}_{st}$ can be obtained from $\tilde{X}_{0,t-s}$ by a time-shift and are therefore also bounded.

By Theorem 3.7 from Ref.~\refcite{franz03b}, $j_{st}(x)=X_{st}=\tilde{X}_{st}\Gamma_{t]}$ for $0\le s\le t$ defines a monotone L\'evy process on $\mathbb{C}[x]$ w.r.t.\ to the vacuum state in the sense of Definition 2.5 in Ref.~\refcite{franz03b}. Since $\tilde{X}_{st}$ and $\Gamma_{t]}$ are symmetric, bounded, and commute for all $0\le s\le t$, the operators $X_t=\tilde{X}_{0t}\Gamma_{t]}$ are also symmetric and bounded, hence self-adjoint. This implies that $(X_t)_{t\ge 0}$ is a monotone L\'evy process in the sense of Definition \ref{mon-def-mon-levy}.

The coefficients $M_x=\pi(x)$, $\mathbf{1}_\mathbb{R}=\eta(x)$, and $-a=L(x)$ in Equation \eqref{mon-bialg-qsde} correspond exactly to the Sch\"urmann triple associated to the characteristic pair $(a,\rho)$ in Lemma \ref{mon-lemma-triple}. Therefore $(X_t)_{t\ge0}$ has the correct monotone convolution semigroup.
\end{proof}

\begin{remark}
It follows from Ref.~\refcite{franz+schuermann+skeide03} that $\Omega$ is cyclic for $(X_t)_{t\ge 0}$. But $\Omega$ is not separating, except for $\rho=0$ (i.e.\ the pure drift process, see Subsection \ref{mon-exa-drift}). Set $\overline{X}_t=X_t-\Phi(X_t)\Gamma_{t]}$ for $t\ge 0$. Then we have
\begin{eqnarray*}
\Phi\left((\overline{X}_s(\overline{X}_t-\overline{X}_s))^*\overline{X}_s(\overline{X}_t-\overline{X}_s)\right) &=& \Phi\left((\overline{X}_t-\overline{X}_s)\overline{X}^2_s(\overline{X}_t-\overline{X}_s)\right) \\
&=& \big(\Phi(\overline{X}_t-\overline{X}_s)\big)^2 \Phi(\overline{X}_s^2)=0
\end{eqnarray*}
and therefore
\[
\overline{X}_s(\overline{X}_t-\overline{X}_s)\Omega =0,
\]
but
\begin{eqnarray*}
\Phi\left(\overline{X}_s\left(\overline{X}_t-\overline{X}_s\right)\left(\overline{X}_s\left(\overline{X}_t-\overline{X}_s\right)\right)^*\right) &=& \Phi\left(\overline{X}_s\left(\overline{X}_t-\overline{X}_s\right)^2\overline{X}_s\right) \\
&=& \Phi\left(\overline{X}_s^2\right)\Phi\left(\left(\overline{X}_t-\overline{X}_s\right)^2\right)
\end{eqnarray*}
proves that $\overline{X}_s(\overline{X}_t-\overline{X}_s)\not=0$ for $0<s<t$, unless $\overline{X}_t=0$ for all $t\ge0$.
\end{remark}

\section{The Markov Semigroup of a Monotone L\'evy Process}

Let $(a,\rho)$, $(a,\rho)\not=(0,0)$, be a non-trivial characteristic pair. In this section we will always assume that $(X_t)_{t\ge0}$ is the monotone L\'evy processes on $\Gamma\Big(L^2\big(\mathbb{R}_+,L^2(\mathbb{R},\rho)\big)\Big)$ constructed in Theorem \ref{mon-theo-mon-levy} for $(a,\rho)$.

As ``conditional expectations'' we will use the linear maps $E_t:\mathcal{B}\Big(\Gamma\big(L^2(\mathbb{R}_+,H)\big)\Big)\to\mathcal{B}\Big(\Gamma\big(L^2(\mathbb{R}_+,H)\big)\Big)$,
\begin{equation}\label{mon-eq-cond-exp}
E_t(X)=\Gamma_{t]}X\Gamma_{t]}, \qquad X\in\mathcal{B}\Big(\Gamma\big(L^2(\mathbb{R}_+,H)\big)\Big),
\end{equation}
for $t\ge0$. Denote the image of $E_t$ by
\[
\mathcal{A}_{t]}=E_t\Big(\mathcal{B}\Big(\Gamma\big(L^2(\mathbb{R}^+,H)\big)\Big)\Big).
\]
These are exactly the operators on the Fock space which are $\Omega$-adapted in the sense of Belton.\cite{belton01,belton03} The algebra $\mathcal{A}_{t]}$ consists of all bounded operators that leave the subspace $\Gamma\big(L^2([0,t],H)\big)\otimes \Omega \subseteq \Gamma\big(L^2([0,t],H)\big)\otimes \Gamma\big(L^2([t,\infty[,H)\big)\cong \Gamma\big(L^2(\mathbb{R}_+,H)\big)$ invariant and vanish on its orthogonal complement.

The conditional expectations have the following properties.
\begin{lemma}
\begin{itemize}
\item[(a)]
$E_s\circ E_t = E_{s} = E_t\circ E_s$, for $0\le s\le t$, and in particular $E_t^2=E_t$.
\item[(b)]
$E_t$ is completely positive for all $t\ge 0$.
\item[(c)]
$E_t(XYZ)=XE_t(Y)Z$ for all $X,Z\in\mathcal{A}_{t]}$, $Y\in\mathcal{B}\Big(\Gamma\big(L^2(\mathbb{R}^+,H)\big)\Big)$.
\item[(d)]
$E_t({\rm id})=\Gamma_{t]}$ for all  $t\ge 0$.
\end{itemize}
\end{lemma}

Let $(k_t)_{t\ge 0}$ be the quantum stochastic process on $\mathbb{C}[x]$ defined by
\[
k_t(x^k)=X_t^k\Gamma_{t]} =\left\{\begin{array}{lcl}
\Gamma_{t]} & \mbox{ if } &k=0, \\
X^k_t & \mbox{ if } & k=1,2,\ldots.
\end{array}\right.
\]
Using functional calculus we extend $(k_t)_{t\ge 0}$ to $C_b(\mathbb{R})$ by $k_t(f)=f(X_t)\Gamma_{t]}$. Denote by $\mathcal{A}_t$ the algebra generated by $\Gamma_{t]}$ and $X_t$, i.e.\ $\mathcal{A}_{t}=k_t(\mathbb{C}[x])$ and by
\[
\hat{A}_t=\{\Gamma_{t]}f(X_t)|f\in C_b(\mathbb{R})\}
\]
the image of the extension of $k_t$ to $C_b(\mathbb{R})$. We have of course
\[
\mathcal{A}_s\subseteq\hat{\mathcal{A}}_s\subseteq\mathcal{A}_{t]}
\]
for all $0\le s\le t$.

\begin{theorem}\label{mon-theo-markovianity}
The monotone L\'evy process is Markovian, i.e.\ we have
\[
E_s\big(k_t(P)\big)\in \mathcal{A}_s \quad \mbox{ and } \quad E_s\big(k_t(f)\big)\in \hat{\mathcal{A}}_s
\]
for all $0\le s\le t$, $P\in\mathbb{C}[x]$, and $f\in C_b(\mathbb{R})$.

The semigroup of Markovian transition operators $T_t:C_b(\mathbb{R})\to C_b(\mathbb{R})$ with
\[
E_s\big(k_t(f)\big)=k_s\big(T_{t-s}(f)\big)
\]
for $0\le s\le t$, $f\in C_b(\mathbb{R})$ is given by
\begin{equation}\label{mon-eq-markov}
T_tf(x)=\int_{\mathbb{R}}f(y){\rm d}\mu_{t,x}(y),
\end{equation}
for $t\ge 0$, $f\in C_b(\mathbb{R})$, where $\mu_{t,x}=\delta_x\triangleright\mu_t$. The semigroup $(T_t)_{t\ge0}$ maps $C_b(\mathbb{R})$ to itself and polynomials to polynomials.
\end{theorem}
\begin{proof}
The computation of $T_{t-s}$ is the same as in the proof of Proposition \ref{mon-prop-cond-exp} and Corollary \ref{mon-cor-cond-exp}, just write $X_t$ as a sum $X_t=X_s+(X_t-X_s)$ of two monotonically independent self-adjoint operators. 

The formula in Proposition \ref{mon-prop-cond-exp} can be interpreted as a generating function and shows that polynomials are mapped to polynomials. On then other hand, equation \eqref{mon-eq-markov} implies that $T_tf$ is again in $C_b(\mathbb{R})$, due to the continuity of the monotone convolution.
\end{proof}

As in Ref.~\refcite{biane98,bozejko+kuemmerer+speicher96}, the Markov property implies the existence of a classical version.
\begin{corollary}
There exists a classical Markov process $(\hat{X}_t)_{t\ge 0}$ on $\mathbb{R}$ that has the same time-ordered joint expectations as $(X_t)_{t\ge0}$, i.e.\
\[
\mathbb{E}\big(f_1(\hat{X}_{t_1})\cdots f_n(\hat{X}_{t_n})\big)=\Phi\big(f_1(X_{t_1})\cdots f_n(X_{t_n})\big)
\]
for all $n\in\mathbb{N}$, $0\le t_1\le \cdots\le t_n$, $f_1,\ldots,f_n\in C_b(\mathbb{R})$.
\end{corollary}

We have the following expression for the generator of the semigroup $(T_t)_{t\ge 0}$.

\begin{proposition}\label{mon-generator of T_t}
Let $(X_t)_{t\ge0}$ be the monotone L\'evy process whose monotone convolution semigroup has characteristic pair $(a,\rho)$. Then the generator
\[
\mathcal{L} = \left.\frac{{\rm d}}{{\rm d}t}\right|_{t=0} T_t
\]
of the associated semigroup of transition operators is given by
\[
\mathcal{L}f(x)=-af'(x) + \int_{\mathbb{R}} \frac{f(x)-f(y)-(x-y)f'(x)}{(x-y)^2}{\rm d}\rho(y)
\]
for $f\in C^2_b(\mathbb{R})$.
\end{proposition}
\begin{proof}
We use $f_z(x)=\frac{1}{z-x}$ as a generating function. We have
\[
T_tf_z(X_s) = E_s\left(\frac{1}{z-X_{s+t}}\right)=\frac{\Gamma_{t]}}{H_{\mu_t}(z)-X_s}=k_s\left(\frac{1}{H_{\mu_t}(z)-x}\right),
\]
and therefore
\[
\mathcal{L}f_z(x)=-\frac{A(z)}{(z-x)^2}= -\frac{a}{(z-x)^2}+\int_\mathbb{R}\frac{1}{(z-x)^2(z-y)}{\rm d}\rho(y).
\]
Using a partial fraction decomposition, this becomes
\begin{eqnarray*}
&& \mathcal{L}f_z(x) = \\
&=&  -\frac{a}{(z-x)^2}+\int_\mathbb{R}\left(\frac{1}{(x-y)^2}\left(\frac{1}{z-x}- \frac{1}{z-y}\right)+\frac{1}{x-y} \frac{1}{(z-x)^2}\right){\rm d}\rho(y) \\
&=& -a f_z'(x) + \int_\mathbb{R}\frac{f_z(x)-f_z(y)+(x-y)f_z'(x)}{(x-y)^2}{\rm d}\rho(y)
\end{eqnarray*}
For sufficiently large $|z|$, the series converge uniformly on a bounded interval containing ${\rm supp}\,\rho$. Therefore we can interchange summation and integration and deduce that the formula given in the proposition holds for polynomials.
Since $\rho$ has compact support, the formula extends to functions in $C_b^2(\mathbb{R})$ by a Stone-Weierstrass type approximation.
\end{proof}

\section{Examples}

\subsection{Pure drift process}\label{mon-exa-drift}

The simplest case is the monotone L\'evy process associated to the characteristic pair $(a,0)$, $a\in\mathbb{R}$. We get $A(z)=a$, $H_t(z)=z-at$, and $\mu_t=-at$. The associated monotone L\'evy process is just $X_t=-a{\rm id}$.

\subsection{Monotone Brownian motion}

Consider now $(0,\delta_0)$. The we get $A(z)=-\frac{1}{z}$, $H_t(z)=\sqrt{z^2-2t}$, and $\mu_t$ is absolutely continuous w.r.t.\ to Lebesgue measure, with density
\[
\frac{1}{\pi\sqrt{2t-x^2}}\mathbf{1}_{]-\sqrt{2t},\sqrt{2t}[}, \qquad t\ge 0.
\]
The generator of the semigroup $(T_t)_{t\ge 0}$ is given by
\[
\mathcal{L}f(x)=\left\{\begin{array}{lcl}
\frac{f(x)-f(0)-xf'(x)}{x^2} & \mbox{ if } & x\not=0, \\
\frac{1}{2}f''(0) & \mbox{ if } & x=0,
\end{array}\right.
\]
on $f\in C_b^2(\mathbb{R})$.

The process $\tilde{X}_{0t}$ is equal to the quantum Az\'ema martingale\cite{parthasarathy90} with parameter $q=0$, cf.\ Remark 4.9 in Ref.~\refcite{franz03b}. The classical version is the classical Az\'ema martingale.\cite{emery89}

The monotone L\'evy process with characteristic pair $(a,\delta_0)$, $a\in\mathbb{R}$, is the monotone analogue of a Brownian motion with drift $-a$. We get $A(z)=a-\frac{1}{z}$. The reciprocal Cauchy transform is given as the unique solution of $H_t(z)=w$ in $\mathbb{C}^+$ of
\[
a(w-z)+\ln\frac{aw-1}{az-1} = a^2t.
\]
The generator of the semigroup $(T_t)_{t\ge 0}$ is given by
\[
\mathcal{L}f(x)=\left\{\begin{array}{lcl}
\frac{f(x)-f(0)-x(1+ax)f'(x)}{x^2} & \mbox{ if } & x\not=0, \\
\frac{1}{2}f''(0)-af'(0) & \mbox{ if } & x=0,
\end{array}\right.
\]
on $f\in C_b^2(\mathbb{R})$.

\subsection{Monotone Poisson process}

Let now $(a,\rho)=\left(-\frac{\lambda}{2},\frac{\lambda}{2}\delta_1\right)$ with $\lambda>0$. Then we have $A(z)=\frac{\lambda}{2}\frac{z}{1-z}$ and $H_t(z)=w$ is the unique solution in $\mathbb{C}^+$ of
\[
-\frac{\lambda}{2}(w-z)-\frac{\lambda}{2}\ln\frac{w-1}{z-1}=t.
\]
The corresponding probability measures where determined in Ref.~\refcite{muraki01}, where they arose as limit distributions in a Poisson-type limit theorem, see also Example 4.4.(3) in Ref.~\refcite{muraki00}. The monotone L\'evy process associated to this characteristic pair is the monotone Poisson process. The generator of its Markov semigroup $(T_t)_{t\ge 0}$ is given by
\[
\mathcal{L}f(x)=\frac{\lambda}{2}\frac{f(x)-f(1)-xf'(x)}{(x-1)^2}
\]
for $f\in C_b^2(\mathbb{R})$.

\section{Martingales}

In this section we show how one can construct a family of martingales from a monotone L\'evy process that is analogous to the family of exponential martingales of a classical L\'evy process. If $(X_t)_{t\ge 0}$ is a classical L\'evy process with characteristic functions $\varphi_t(u)=\mathbb{E}(e^{iuX_t})$, then for any $u\in\mathbb{R}$, the process $(M^u_t)_{t\ge 0}$ with
\[
M^u_t=\exp\big(iuX_t-\log \varphi_t(u)\big),\qquad \mbox{ for } t\ge 0,
\]
is a martingale w.r.t.\ the filtration of $(X_t)_{t\ge 0}$. An analogous family of martingales for free increment processes has been defined by Biane, see Section 4.3 in Ref.~\refcite{biane98}.

\begin{definition}
Let $H$ be a Hilbert space, $\Gamma\big(L^2(\mathbb{R}_+,H)\big)$ the Fock space over $L^2(\mathbb{R}_+,H)$, and $(E_t)_{t\ge 0}$ the family of conditional expectations introduced in \eqref{mon-eq-cond-exp}.

We call a family $(M_t)_{t\in I}$ of operators on $\Gamma\big(L^2(\mathbb{R}_+,H)\big)$ indexed by an interval $I\subseteq\mathbb{R}_+$ a {\em martingale}, if
\[
E_s(M_t)=M_s
\]
holds for all $s,t\in I$ with $s\le t$.
\end{definition}

\begin{lemma}\label{mon-lem-inj}
Let $(\mu_t)_{t\ge 0}$ be a monotone convolution semigroup of probability measures with reciprocal Cauchy transforms $(H_t)_{t\ge 0}$. Then the $H_t$ are injective on $\mathbb{C}^+$ for all $t\ge 0$.
\end{lemma}
\begin{proof}
If $\mu_t=\delta_0$ for all $t\ge 0$, then $H_t={\rm id}$ for all $t\ge 0$ and the lemma is true.

Assume now that $(\mu_t)$ is a non-trivial monotone convolution semigroup.

Let $t\ge 0$ and $z_1,z_2\in H_t(\mathbb{C}^+)$ with $H_t(z_1)=H_t(z_2)$. Then we have
\[
\int_{z_1}^w \frac{{\rm d} z}{A(z)} = t = \int_{z_2}^w\frac{{\rm d}z}{A(z)},
\]
where $w=H_t(z_i)$, $i=1,2$, and $A(z)$ denotes the generator of $(H_t)_{t\ge 0}$, see Theorem 4.7 in Ref.~\refcite{muraki00} (or also Theorem \ref{mon-theo-levy-khintchine} for the case where the $\mu_t$ are compactly supported). This implies
\[
\int_{z_1}^{z_2} \frac{{\rm d} z}{A(z)}=0,
\]
and $z_1=H_0(z_1)=z_2$ by the uniqueness of the solution in Theorem 4.7.(2) in Ref.~\refcite{muraki00} or Theorem \ref{mon-theo-levy-khintchine}(b).
\end{proof}

\begin{theorem}\label{mon-theo-martingale}
Let $T>0$, and let $(\mu_t)_{t\ge 0}$ be a monotone convolution semigroup of compactly supported probability measures with reciprocal Cauchy transforms $(H_t)_{t\ge 0}$ and monotone L\'evy process $(X_t)_{t\ge 0}$. Then for any $z\in H_T(\mathbb{C}^+)$, the operator process $(M^z_t)_{0\le t\le T}$ with
\[
M^z_t=k_t\left(\frac{1}{H_t^{-1}(z)-X_t}\right) =\frac{1}{H_t^{-1}(z)-X_t}\Gamma_{t]}
\]
is a martingale.
\end{theorem}
\begin{proof}
$M_t^z$ is well-defined, since we have $H_T(\mathbb{C}^+)\subseteq H_t(\mathbb{C}^+)$ for $0\le t\le T$ by the semigroup property of $(H_t)_{t\ge 0}$ and since the $H_t$ are injective by Lemma \ref{mon-lem-inj}.

Let now $0\le s\le t\le T$, then we have
\begin{eqnarray*}
E_s(M^z_t) &=& E_s \left(\frac{1}{H_t^{-1}(z)-X_t}\Gamma_{t]}\right) = \Gamma_{s]}\frac{1}{H_t^{-1}(z)-X_s-(X_t+X_s)}\Gamma_{s]} \\
&=& \frac{1}{H_{t-s}\big(H_t^{-1}(z)\big)-X_s}\Gamma_{s]}= \frac{1}{H_s^{-1}(z)-X_s}\Gamma_{s]} \\
&=& M^z_s,
\end{eqnarray*}
cf.\ Proposition \ref{mon-prop-cond-exp}.
\end{proof}

\end{document}